\documentclass[conference]{IEEEtran}
\IEEEoverridecommandlockouts
\usepackage{cite}
\usepackage{amsmath,amssymb,amsfonts}
\usepackage{algorithmic}
\usepackage{graphicx}
\usepackage{textcomp}
\usepackage{xcolor}
\def\BibTeX{{\rm B\kern-.05em{\sc i\kern-.025em b}\kern-.08em
    T\kern-.1667em\lower.7ex\hbox{E}\kern-.125emX}}

\newtheorem{assumption}{\bf{Assumption}}

\newcommand{\xddots}{%
  \raise 4pt \hbox {.}
  \mkern 6mu
  \raise 1pt \hbox {.}
  \mkern 6mu
  \raise -2pt \hbox {.}
}

\begin{document}

\title{Fault-Ride-Through (FRT) Control of a grid-connected Fixed-Speed Wind Energy Conversion System using STATCOM \\
}
\author{\IEEEauthorblockN{
\textbf{Ayobami Olajube, Satish Vedula, Koto Omiloli, Olugbenga Moses Anubi}
}\\
\IEEEauthorblockA{Department of Electrical and Computer Engineering, FAMU-FSU College of Engineering \\
Center for Advanced Power Systems, Florida State University \\
E-mail: \{aao21b, svedula, kao23a, oanubi\}@fsu.edu}
}
\maketitle

\begin{abstract}
Wind energy conversion system (WECS) is stochastic in nature and has low inertia leading to grid voltage instability, poor reactive power compensation and most importantly fault susceptibility. Variable speed WECS such as the doubly-fed induction generators (DFIG) are well known to reach steady state quickly after fault occurrence without the need for an external reactive power source   because of the presence of a back-to-back converter that provides independent control of the active and reactive power unlike in the fixed-speed squirrel cage induction generator (SCIG) counterpart that can't be stabilized unless an external source of reactive power support is present. However, controlling DFIG is complicated and costly due to complete tripping unlike the fixed-speed generators which doesn't trip completely when fault occurs. Hence, in this work, a 48-pulse, 3-phase static synchronous compensator (STATCOM) is used to ensure reactive power compensation and fault-ride through (FRT) control of the SCIG against over-voltage emanating from fault occurrence in a grid-connected power system. The goal here is to guarantee voltage stability and fault-ride-through (FRT) control against injected faults within certain time ranges at the point of common coupling (PCC) between the AC source, the load, and the fixed-speed WECS. The numerical simulation shows the fault-ride-through capability of the STATCOM-controlled fixed-speed WECS, good voltage regulation and effective reactive power compensation.
\end{abstract}


\section{Introduction}
\begin{figure}[t!]
      \centering
      \includegraphics[width=0.42\textwidth]{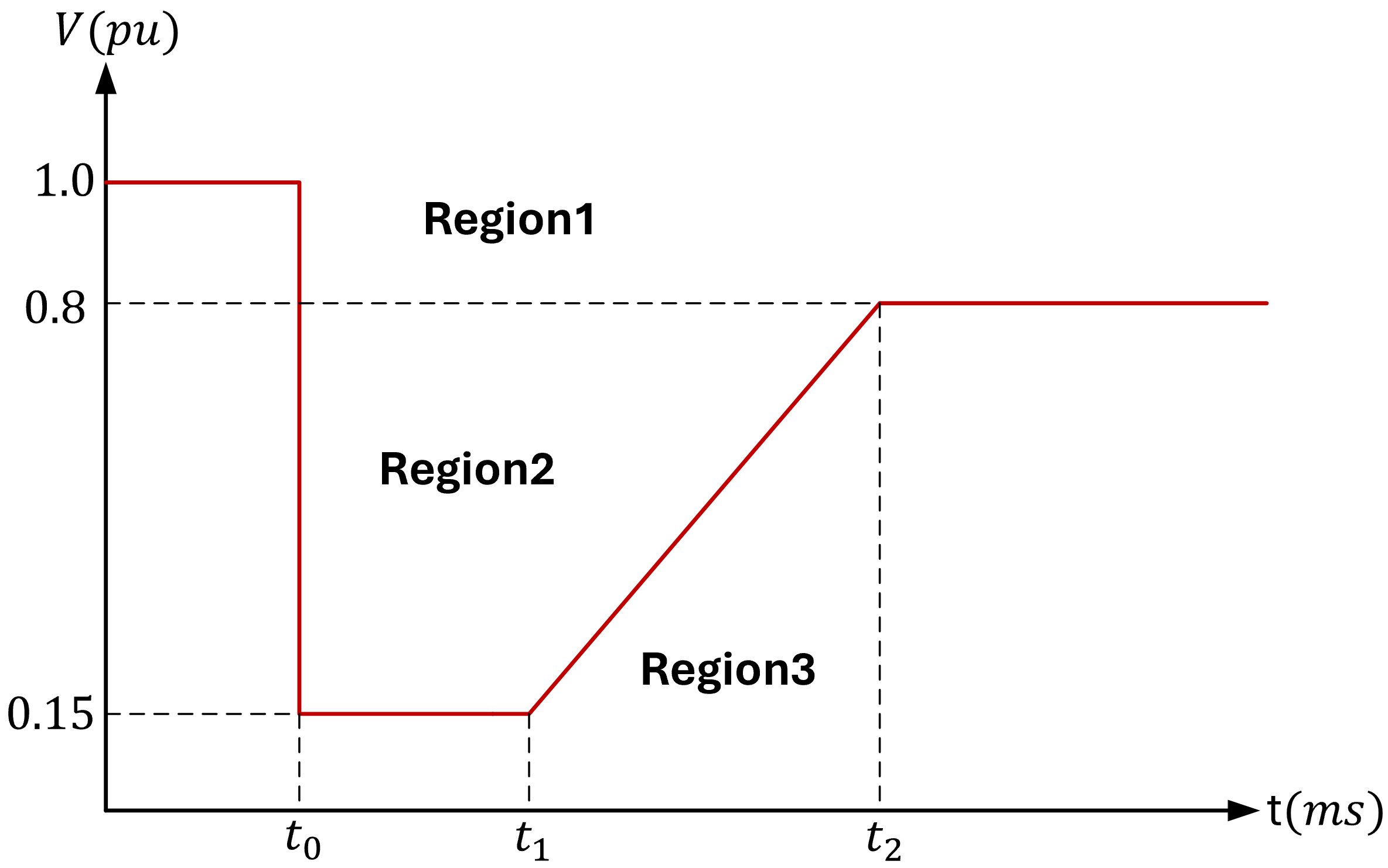}
	 \caption{Generic WECS FRT curve}
     \label{lvrt} 
\end{figure}

{M}{odern} power grid network is gradually undergoing significant structural modifications in other to achieve net-zero power generation. Efforts are in top-gear to achieve almost 100\% clean energy generation and the conventional synchronous generators are gradually being replaced by the low-inertia inverter-based resources to guarantee cleaner energy supply and improve environmental conditions. However, replacing the high-inertia synchronous generators with low-inertia and fast-switching inverter-tied renewable sources will come at cost on voltage and frequency stability, reactive power compensation, and grid robustness against disturbances or sudden load changes \cite{milano2018foundations}. The most popular among the renewable energy sources in a global space is the wind energy conversion system especially in Europe and the united states. The world has witnessed tremendous wind energy penetration in the last few decades with China leading the race in terms of cummulative wind power capacity of 236GW, followed closely by the United States with 106GW capacity~\cite{wiser2020wind} as at the end of 2019 and, this trend will continue to increase. Most wind energy conversion systems such as the squirrel cage induction generators (SCIGs) exist as fixed speed generators. They have simple structure, cheap, rugged and devoid of maintenance requirements in most cases but suffer from stability issues, which necessitates the need for the use of flexible alternating current transmission system (FACTS) devices in order to meet the grid code requirement during faults\cite{operator2011wind,rahimi2010grid}. A generic grid code curve for WECS is shown in Figure~\ref{lvrt} where the wind energy system operates normally at voltage $1.0\mathrm{~p.u}$ at the point-of-common coupling (PCC) in region 1. When a fault occurs at $t_0$, the WECS experiences voltage sag down to $0.15\mathrm{~pu}$ at the PCC and remains connected to the grid until $t_1$ within region 2. However, if the WECS undergoes a further decrease in voltage, it will disconnect completely from the grid due to the lack of reactive power support. Nonetheless, the fixed speed WECS is expected to experience an increase in voltage support in the presence of FACTS devices for restoration to about $0.8 \mathrm{~pu}$ in region~3 \cite{tarafdar2019review}. As the wind turbine penetration increases globally, different countries revise grid codes to suit their grid specifications and requirements \cite{tsili2009review,taul2019current}.   

SCIGs connected directly to stiff grid operate at fixed speeds and require reactive power support. They consume reactive power and slow down voltage restoration at weak nodes during faults, leading to rotor angle and voltage instability \cite{rathi2005novel}. When a fault occurs, the wind turbine generator will draw more current and accelerate faster due to the imbalance between the mechanical power extracted from the wind and the electrical power delivered to the grid. When the fault is cleared, the generator slows down voltage recovery by consuming a large amount of reactive power leading to instability, especially at weak grid level. Hence, fast-switching FACTS devices such as static compensators (STATCOMs) are required to be connected to the system to provide support for voltage restoration \cite{el2008evaluation}. STATCOM is a popular shunt-connected reactive power compensating FACTS device that can provide both transient and dynamic stability in weak grids during fault occurrences \cite{saad1998application, muyeen2005stabilization}. It also helps in power quality improvement in modern grids \cite{kuiava2009control,teleke2011application,barrado2010power}. The variability and the low energy density of the DC-link capacitor of the popular voltage-source-converter (VSC) based STATCOM is its main drawback and energy storage technologies are often added to optimize its performance \cite{hingorani1999understanding}. Other FACTS devices such as thyristor-controlled series capacitor (TCSC), static synchronous series compensator (SSSC), static VAR compensator (SVC), and universal power flow controller (UPFC) have also been used in conjunction with WECS at the PCC to provide both voltage and rotor speed stability depending on the functional requirement, mode of connection and control objective \cite{adetokun2021application}. Control methods such as PI \cite{ledesma2005doubly} , nonlinear \cite{nemmour2010advanced}, robust \cite{sakamoto2007output} and adaptive \cite{jurado2003adaptive} have been used for fault-ride-through (FRT) analysis of the conventional VSC-based STATCOM in fixed-speed WECS but harmonics are reportedly high because of the poor switching dynamics of the VSCs. Hence, the advent of 48 pulse statcom switching technology has helped with keeping the DC-link voltage constant without the need for energy storage system, maintaining voltage balancing and eliminating harmonics \cite{singh2012new,abdellaoui2015statcom}.

The main contribution of this paper is the development of a fault-ride-through (FRT) capability for the grid-connected WECS using 48 pulse STATCOM via cascaded PI-based outer loop DC voltage and the inner loop current controllers. The 48-pulse based STATCOM provides reactive power compensation to the induction generator during fault and also eliminates harmonic distortion that may result from the cascaded control structure and the inverter PWM technique, thereby keeping the voltage profile at the nominal level.

The rest of the paper is organized as follows: section \ref{sec:notations} shows the preliminaries and notations used throughout this paper, section \ref{sec:system} gives general system description, section \ref{sec:models} shows the model development of a wind turbine, induction generator and the STATCOM, section \ref{sec:lvrt control} depicts the FRT control design and algorithm development, section \ref{sec:simulation results} enumerates the numerical simulation results and the conclusion follows in section \ref{sec: conclusion}.

\section{NOTATIONS AND PRELIMINARIES}\label{sec:notations}

In this paper, we use $\mathbb{R}$, $\mathbb{R}^n$, $\mathbb{R}^{n\times m}$ to denote the space of real numbers, real vectors of length $n$, and real matrices of $n$ rows, and $m$ columns, respectively.
$\mathbb{R}_+$ indicates a set of positive real numbers.
$X^\top$ denotes the transpose of the quantity $X$.
Normal-face lower-case letters ($x\in\mathbb{R}$) are used to represent real scalars, bold-face lower-case letters ($\mathbf{x}\in\mathbb{R}^n$) represent vectors, while normal-face upper case ($X\in\mathbb{R}^{n\times m}$) represents matrices. $X\succ 0 \hspace{1mm} (\succeq0)$ denotes a positive definite (semi-definite) matrix. Furthermore, $J=\begin{bmatrix}
    0&1\\
    -1&0\\
    \end{bmatrix}$ is a skew-symmetric matrix,  $\mathbf{1}_n$ is a column vector of ones, and the identity matrix of size $n$ is denoted by $I_n$.
   
To transform variables from three-phase $x_{abc}$ coordinate to $x_{dq0}$ synchronous reference coordinate, we use (\ref{park}) proposed by H. Park in 1929 \cite{park1929two}.
The dq0 transformation matrix for an angle $\theta(t)$ given by  
\begin{equation} \label{park}
   T(\theta)=\sqrt{\frac{2}{3}}\left[\begin{array}{ccc}
\sin \theta & \sin \left(\theta-\frac{2 \pi}{3}\right) & \sin \left(\theta+\frac{2 \pi}{3}\right) \\
\cos \theta & \cos \left(\theta-\frac{2 \pi}{3}\right) & \cos \left(\theta+\frac{2 \pi}{3}\right) \\
\frac{1}{\sqrt{2}} & \frac{1}{\sqrt{2}} & \frac{1}{\sqrt{2}}
\end{array}\right] 
\end{equation}
and satisfies the property $T^{\top}(\theta)=T^{-1}(\theta)$.
Symmetric AC three-phase signals  $x_{a b c}: \mathbb{R} \rightarrow \mathbb{R}^3$ denotes a symmetric three-phase AC signal of the form
\begin{equation}
    x_{a b c}(t)=X_m(t)\left[\begin{array}{c}
\sin \theta \\
\sin \left(\theta-\frac{2 \pi}{3}\right) \\
\sin \left(\theta+\frac{2 \pi}{3}\right)
\end{array}\right]
\end{equation}
where $X_m: \mathbb{R} \rightarrow \mathbb{R}$ is the amplitude and the angle $\theta: \mathbb{R} \rightarrow \mathbb{R}$ satisfies
\begin{equation}
    \dot{\theta}=\omega
\end{equation}
where $\omega: \mathbb{R}\rightarrow \mathbb{R}$ is the grid frequency.

Hence,
\begin{equation}
   x_{d q 0}(t) =   T(\theta)  x_{a b c}(t)
\end{equation}

\begin{figure}[t!] 
	\centering
	\includegraphics[width=0.45\textwidth]{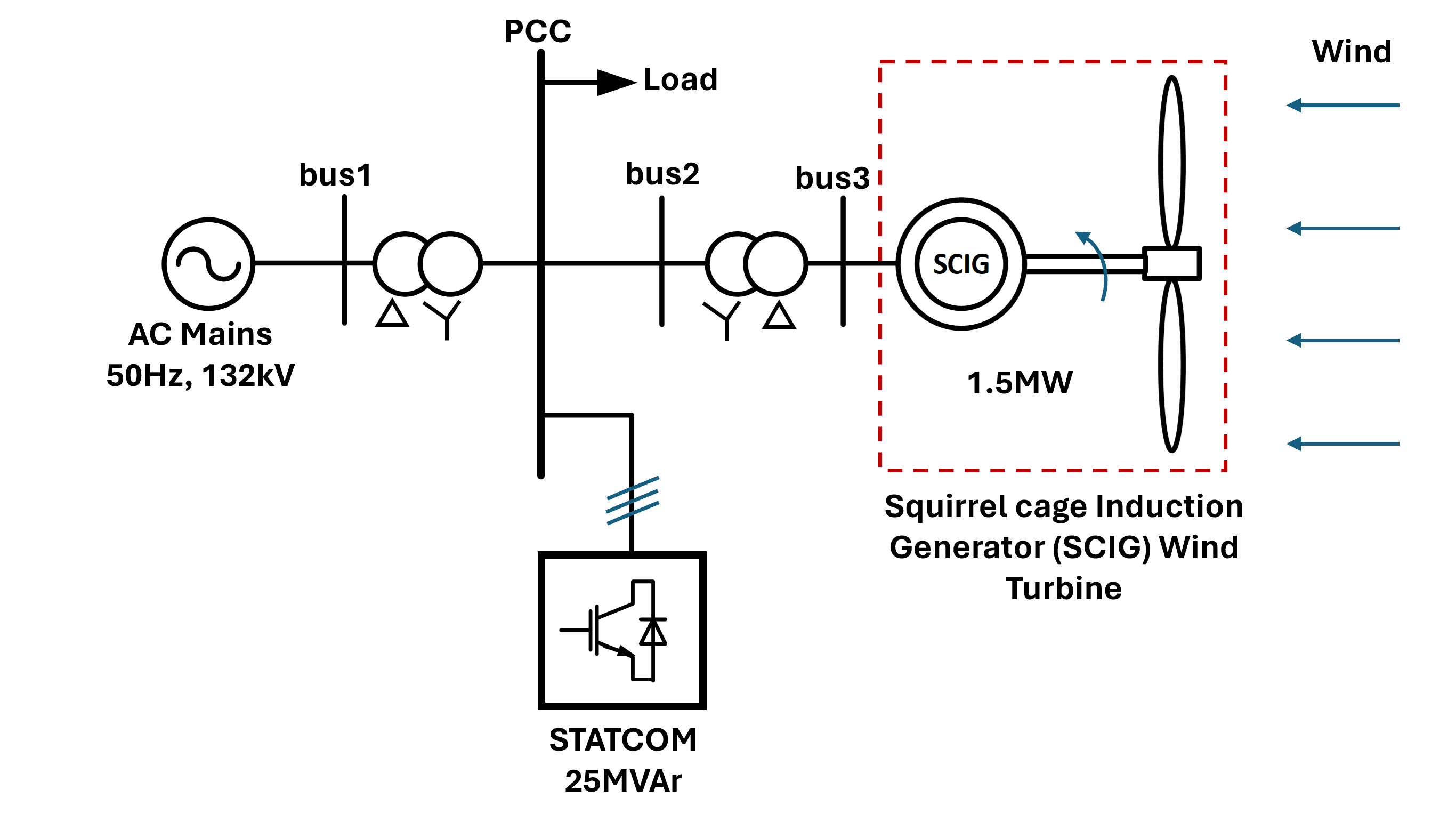} 
	\caption{System Description 
    }
	\label{system}
\end{figure}

 \section{System Description}\label{sec:system}
In this paper, a 1.5~MW SCIG is connected to a medium voltage power system through a transformer. A shunt-connected 25~MVAr STATCOM and a 3-phase constant power load are connected at the PCC and the parameters of the entire power system shown in Figure~\ref{system} are displayed in Table 1.

\begin{table}[ht]
\centering
\caption{Rated Values and System Parameters}
\label{tab:rated}
\resizebox{0.9\columnwidth}{!}{\begin{tabular}{c|c|c}
\hline \hline
\textbf{Parameter}  & \textbf{Parameter}  & \textbf{Parameter}  \\
\textbf{Description} & \textbf{Notation}   & \textbf{Value}  \\  \hline \hline
Wind Turbine Capacity  & $P_{m,rated}$   & 1.5 MW \\ 
Rated wind speed   & $\nu_{rated}$   & $12~ms^{-1}$ \\ 
Inertia Constant & $H$ & $64~sW/VA$ \\
Voltage rating of SCIG & $v$ & $500~V$ \\
Grid frequency & $f$ & $50~Hz$ \\ 
Stator resistance & $r_s$ & $0.01~\Omega$ \\
Rotor resistance & $r_r$ & $0.01~\Omega$ \\
Stator inductance & $l_s$ & $0.041~H$ \\
Rotor inductance & $l_r$ & $0.041~H$ \\
Mutual inductance & $l_m$ & $0.035~H$ \\
Power factor & $p.f$ & $0.85$ \\
Transmission line impedance & $Z_t$ & $(0.12-j2.78)~\Omega/km$ \\ 
Three-phase Constant Load & $S_L$   & $(80+j10)~MVA$ \\ \hline
\end{tabular}}
\end{table}

\section{Model Development}\label{sec:models}
The model of the wind turbine, the SCIG, and the STATCOM are developed in this section:
\subsection{Wind Turbine Model}
\begin{figure}[t!]
      \centering
      \includegraphics[width=0.40\textwidth]{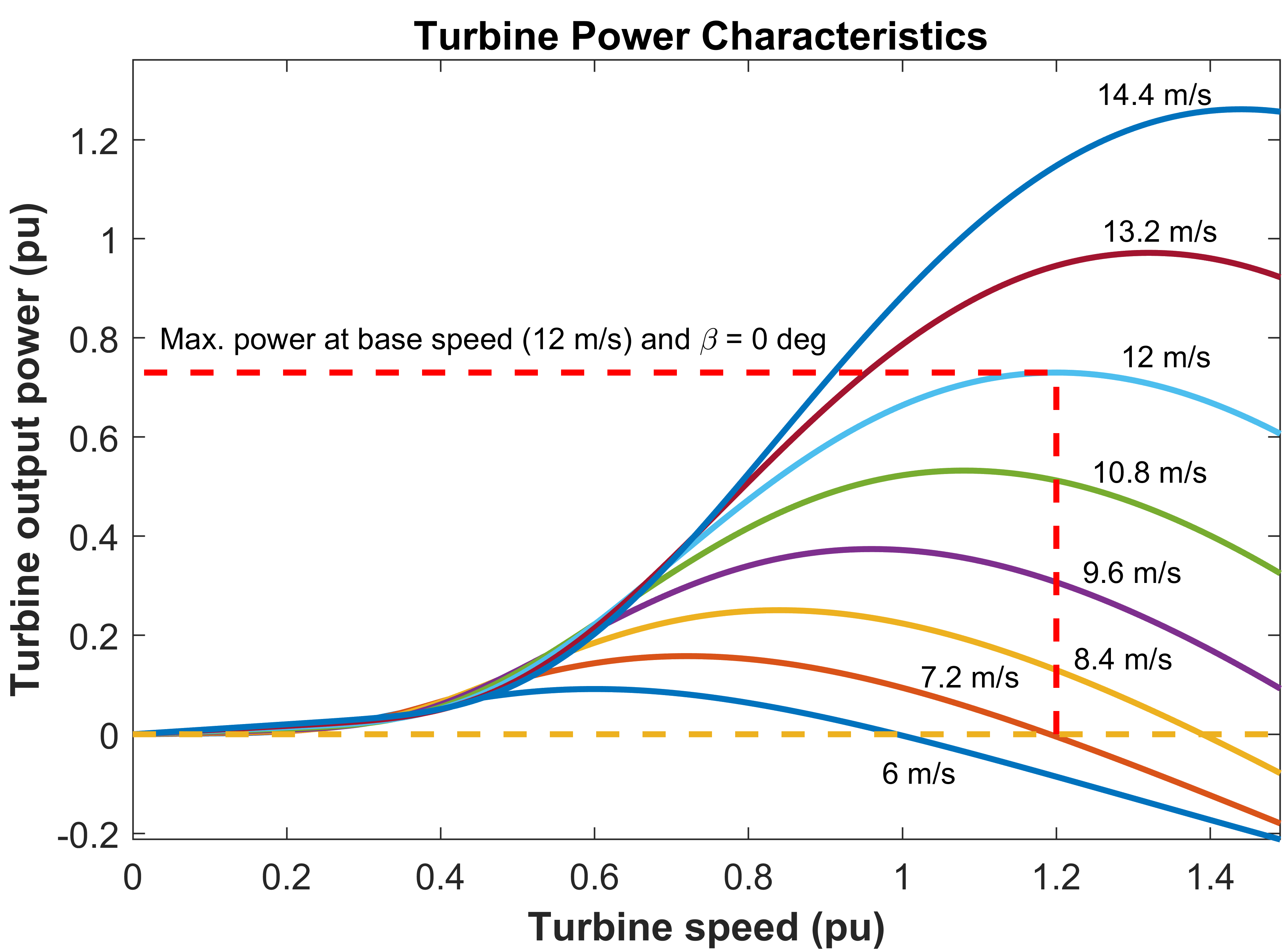}
	 \caption{Turbine Mechanical Output and Speed Characteristics}
     \label{Cp_beta} 
\end{figure}

Based on the aerodynamic model of a wind turbine, the mechanical power that can be extracted from wind by the wind turbine is given as \cite{miller2003dynamic}
\begin{equation}
P_m=\frac{1}{2} \rho A_w \nu^3 C_P(\lambda, \beta)
\end{equation}
where $\rho$ is the air density in $kgm^{-3}$, $A_w$ is the area swept by the turbine rotor blades in $m^2$, $\nu$ is the wind speed in $ms^{-1}$, and $C_P(\lambda, \beta)$ is the power coefficient, which is a function of both the rotor blade pitch angle $\beta$ and the tip-speed ratio $\lambda$. $C_P$ is maximum at a wind speed of $12ms^{-1}$ as shown in Figure~\ref{Cp_beta}. Tip-speed ratio is defined as
\begin{equation}
\lambda=\omega_T L / \nu
\end{equation}
where $\omega_T$ is the wind turbine rotor speed in $rads^{-1}$, $L$ is rotor blade length in meters ($m$), and $\nu$ is the wind speed in $ms^{-1}$. Nonetheless, our attention in this work is focused on the electrical model of the WECS but more details on the control of the wind turbine and pitch angle from mechanical perspective can be studied in \cite{ameli2022hierarchical}.

\subsection{Squirrel Cage Induction Generator Model (SCIG)}
The mathematical model of SCIG with a shorted rotor in dq0 rotational reference coordinate is given as follows \cite{10252896}:
\begin{equation}
\dot{\mathbf{\Phi}} = F\mathbf{\Phi} - N\mathbf{i}+\mathbf{v}
\end{equation}
\begin{equation} \label{swing equation}
\frac{d \omega_r}{\mathrm{dt}}=\frac{p}{ J}\left(T_{\mathrm{e}}-T_m\right)
\end{equation}
\begin{equation} \label{torque equation}
T_{e}=\frac{3 p}{4} l_m\left(\mathbf{i}_{d q s}^T J \mathbf{i}_{d q r}\right)
\end{equation}

\begin{figure*}[t!]
      \centering
      \includegraphics[width=0.65\textwidth]{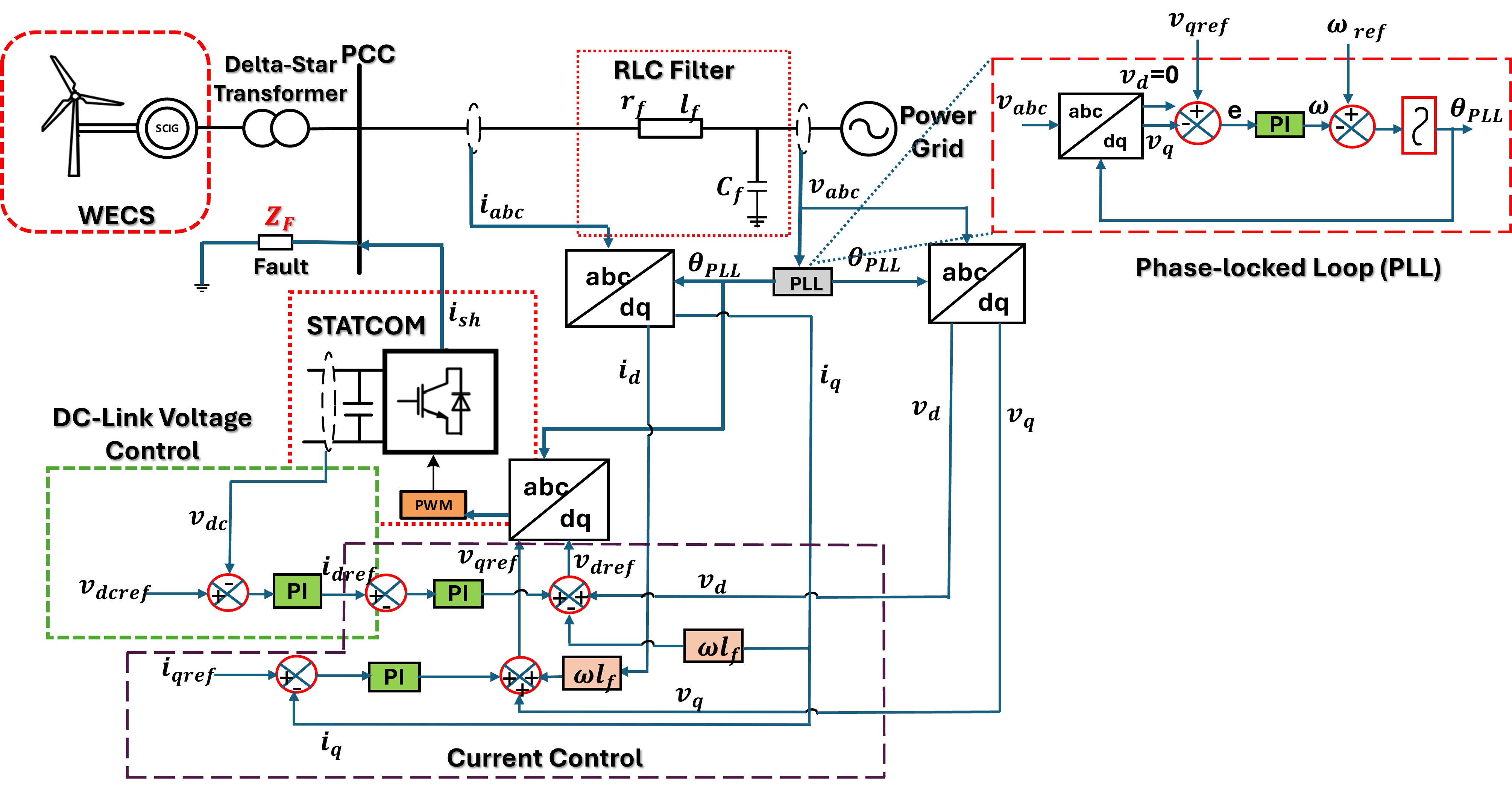}
	 \caption{FRT Control Diagram}
     \label{lvrt control} 
\end{figure*}

where $\mathbf{\Phi} = \left[\mathbf{\psi}_{ds}, \mathbf{\psi}_{qs},\mathbf{\psi}_{dr}, \mathbf{\psi}_{qr}\right]^{\top}$ is the vector of stator and the rotor fluxes for both the direct and quadrature axes in dq0 coordinate, $\mathbf{i} = \left[i_{ds}, i_{qs},i_{dr}, i_{qr}\right]^{\top}$ is the vector of stator and rotor currents, $\mathbf{v} = \left[v_{ds}, v_{qs},0, 0\right]^{\top}$ is the vector of stator and rotor voltages with $v_{dr}=v_{qr}= 0$ because of the shorted rotor condition in SCIG, $N =\textsf{diag}\left(r_s ,r_s ,r_r ,r_r \right)$ is the diagonal matrix of stator and rotor resistances, and

\begin{equation*}
F=\left[\begin{array}{cccc}
0 & \omega & 0 & 0 \\
-\omega & 0 & 0 & 0 \\
0 & (\omega-\omega_{r}) & 0 & 0 \\
-(\omega-\omega_{r})  & 0 & 0 & 0
\end{array}\right]
\end{equation*}
with $\omega, \omega_{r}\in \mathbb{R}_+$ representing the stator and rotor angular frequencies in $\mathrm{~rads^{-1}}$ respectively. The swing equation is represented in (\ref{swing equation}), where  $\frac{d \omega_r}{\mathrm{dt}}$ is the rate of change of rotor frequency in$\mathrm{~Hzs^{-1}}$, $J$ is the generator inertia in$\mathrm{~kgm^2}$, $\rho$ is the number of poles of the generator, $T_e$ is the electromagnetic torque and $T_m$ is mechanical torque from the wind turbine, both in$\mathrm{~Nm}$. $l_m$ from (\ref{torque equation}) is the mutual inductance, $\mathbf{i}_{d q s}^T = \left[i_{ds}, i_{qs}\right]$ and $\mathbf{i}_{d q r} = \left[i_{dr}, i_{qr}\right]^{\top}$. 

The relationship between the fluxes and the currents is given by
\begin{equation*}
    \left\lbrack \begin{array}{c}
\psi_{ds} \\
\psi_{qs} \\
\psi_{dr} \\
\psi_{qr} 
\end{array}\right\rbrack =\left\lbrack \begin{array}{cccc}
l_s  & 0 & l_m  & 0\\
0 & l_s  & 0 & l_m \\
l_m  & 0 & l_r  & 0\\
0 & l_m  & 0 & l_r 
\end{array}\right\rbrack\left\lbrack \begin{array}{c}
i_{ds} \\
i_{qs} \\
i_{dr} \\
i_{qr} 
\end{array}\right\rbrack
\end{equation*}
where $l_s, l_r, l_m$ are the stator, rotor, and mutual inductances of the SCIG respectively.

\begin{figure*}[h!]
      \centering \includegraphics[width=0.67\textwidth]{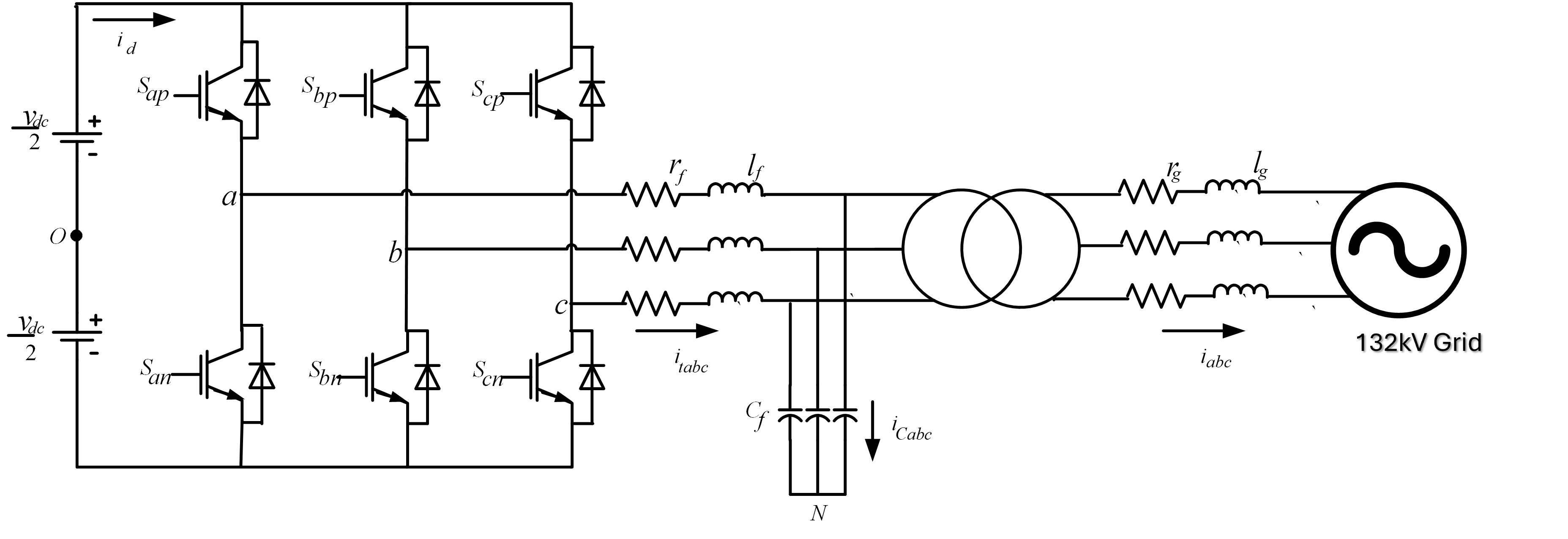}
	 \caption{STATCOM Circuit Diagram}
     \label{statcom} 
\end{figure*}

\begin{figure}[h!]
      \centering \includegraphics[width=0.40\textwidth]{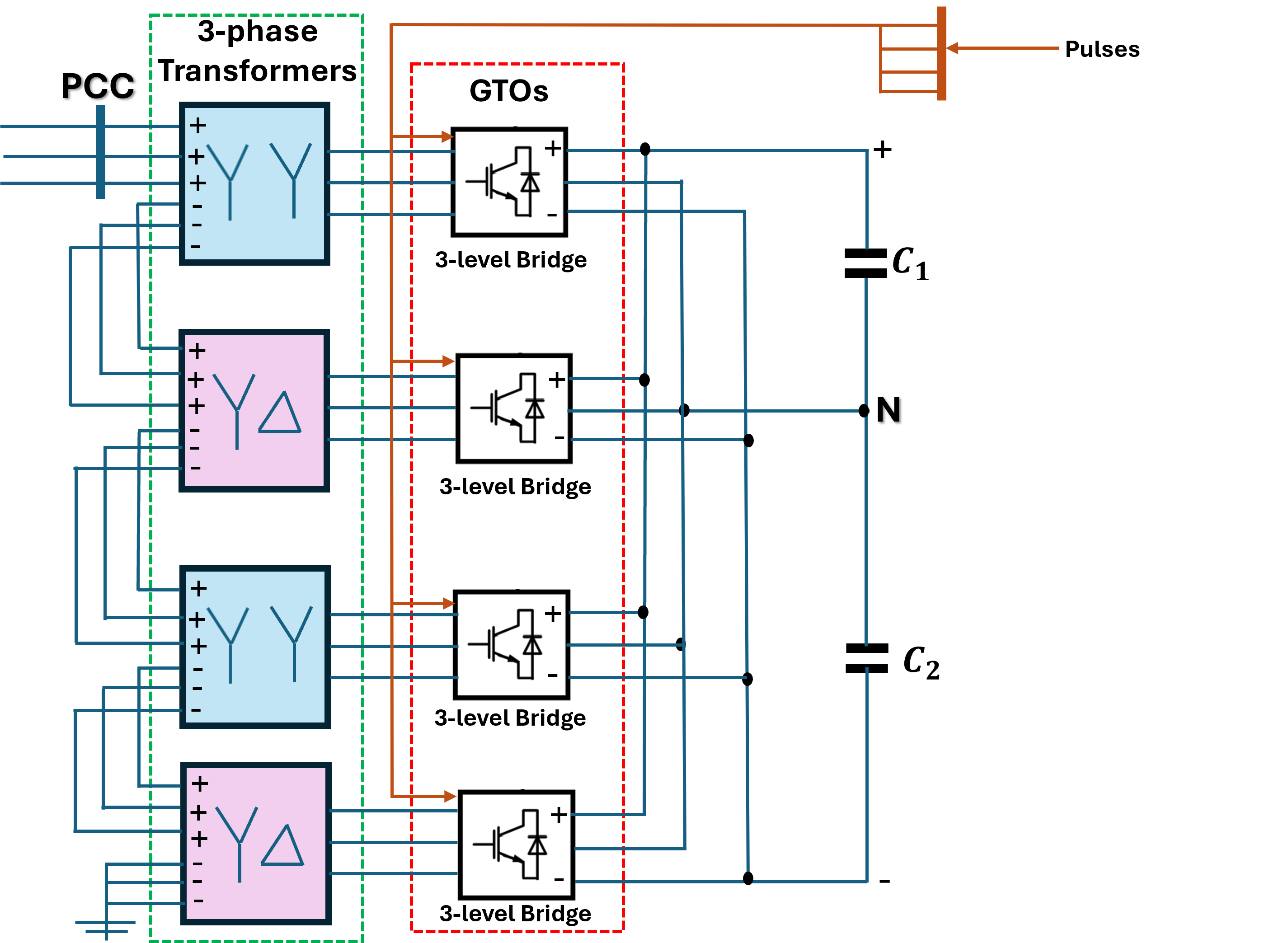}
	 \caption{48 Pulse STATCOM with Phase Shifting Transformers}
     \label{48Pulsestatcom} 
\end{figure}

\subsection{STATCOM Model}
The state space dynamics of a STATCOM shown in Figure~\ref{statcom} is transformed to dq0 coordinate from abc coordinate using park transformation matrix in (\ref{park}) as follows; 
\begin{equation} \label{statcom dynamics}
l_f \frac{d \mathbf{i_t}}{\mathrm{dt}}=-\left(r_f-\omega l_f J\right) \mathbf{i_t}+\mathbf{v_t}-\mathbf{v}
\end{equation}
where $l_f$ is the filter inductance in$\mathrm{~H}$, $r_f$ is the filter resistance in$\mathrm{~\Omega}$, $\omega$ is the grid frequency in $\mathrm{~rads^{-1}}$, which is the same as the stator frequency of the generator, $\mathbf{i_t} = \left[i_{td}, i_{tq}\right]^{\top}$ is the vector of filter currents in dq0 coordinate, $\mathbf{v_t} = \left[\mathbf{v_{td}}, \mathbf{v_{tq}}\right]^{\top}$ is the vector of inverter output voltages in dq0 coordinate, and $\mathbf{v} = \left[\mathbf{v_{d}}, \mathbf{v_{q}}\right]^{\top}$ is the vector of grid voltages in dq0 coordinate. The STATCOM used in this work is a 48 pulse 25~MVAr configuration and the advantage of the 48-pulse STATCOM shown in Figure~\ref{48Pulsestatcom} is the ease and accuracy of generating sinusoidal waves without harmonic distortion or power quality issues \cite{sahoo2006modeling}. The circuit use gate-turn-off (GTO) switches and phase-shifting transformers connected in zigzag form to achieve voltage regulation, fault-ride-through and reactive power compensation. 
\begin{assumption}
  \it{We simplify equation~(\ref{statcom dynamics}) with the assumption of a balanced system and the fact that the filter capacitor $C_f$ is taken to be very small. Hence, the filter capacitor voltage is approximately equal to the grid voltage $\mathbf{v}$.}
\end{assumption}

\section{FRT Control}\label{sec:lvrt control}
Figure~\ref{lvrt control} shows the complete control diagram to enable the fault-ride-through capability of a STATCOM-controlled grid-tied WECS. The main components of this control architecture are the DC-link voltage controller, the current controller, and the phase-locked loop (PLL). The PLL generates the grid phase angle for the coordinate transformation from abc to dq. The DC-link voltage controller keeps the DC link voltage regulated to the constant DC reference voltage value specified at the STATCOM DC input as $v_{dcref}$ using PI controller. The output current $i_{dref}$ of the DC link controller is utilized as a reference for the inner loop current controller. The comparison of the reference $v_{ref}$ from PLL circuit and that of the grid or filter voltage $v_{g}$ are passed through a PI controller to generate reactive current reference $i_{qref}$. The comparison between the current reference $i_{qref}$ and the injected reactive current $i_{sh}$ will produce the phase-shift of the STATCOM voltage from the grid voltage \cite{suresh2010dynamic}. Also, the comparison of the current references via a PI controller is used to generate voltage references $v_{dref}$ and $v_{qref}$, which are transformed into 3-phase voltages $v_{abcref}$ which are compared with the filter output voltages  $v_{abc}$ to generate PWM signal that will switch the STATCOM. The STATCOM in turn injects or absorbs reactive power to maintain the grid voltage profile constant during fault occurrence.
The gains are appropriately selected for the PI controllers in both the outer DC voltage regulation loop and the current control loop respectively as shown in (\ref{idref}) and (\ref{iqref}).

\begin{equation} \label{idref}
    i_{dref} =  K_{p}(v_{dcref} -v_{dc})+ K_{i} \int_{0}^{t} (v_{dcref}(\tau) -v_{dc}(\tau))\,d{\tau} \
\end{equation}

\begin{equation}\label{iqref}
    i_{qref} =  K_{p}(v_{ref} -v_{g})+ K_{i} \int_{0}^{t} (v_{ref}(\tau) -v_{g}(\tau))\,d{\tau} \
\end{equation}

\section{Numerical Simulation}\label{sec:simulation results}
\begin{figure}[t!]
      \centering      \includegraphics[width=0.42\textwidth]{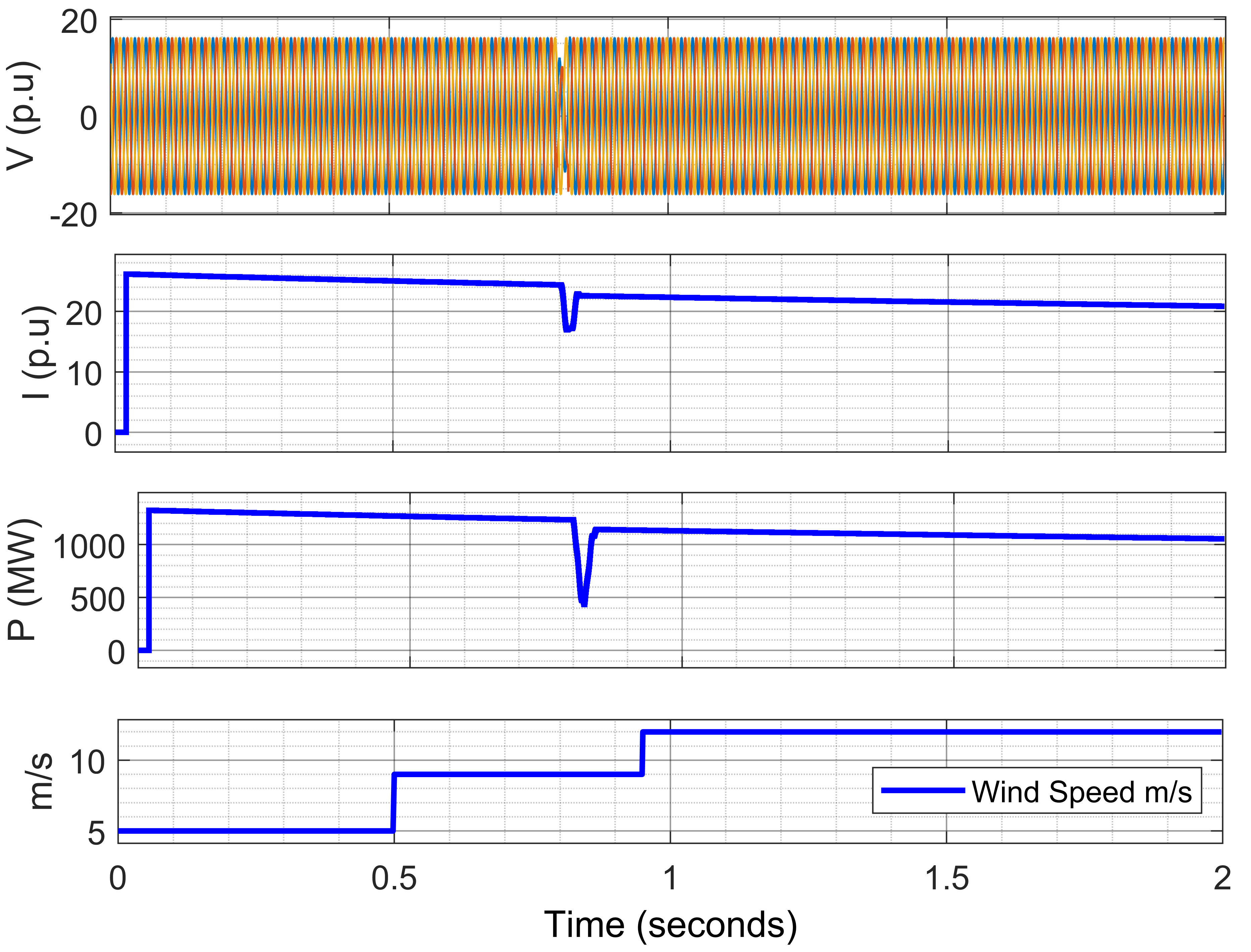}
	 \caption{WECS response without STATCOM}
     \label{nostat_wecs} 
\end{figure}

Using the parameters shown in Table~\ref{tab:rated}, we carried out numerical simulations of the grid-connected wind energy conversion system to obtain the results shown in this section. Initially, without STATCOM being connected, a three-phase fault is applied at the PCC between $0.8\mathrm{~secs}$ and $0.82\mathrm{~secs}$. Figure~\ref{nostat_wecs} shows the variation in the wind speed applied to the system with step increases from $3\mathrm{~ms^{-1}}$ and remain constant between zero to $0.5\mathrm{~secs}$ reaching  $12\mathrm{~ms^{-1}}$ between $1\mathrm{~secs}$ to $2\mathrm{~secs}$. There is no voltage regulation at the terminal of the WECS, because it has over 1000\% overvoltage when compared to the expected nominal value of $1\mathrm{~p.u}$. The WECS becomes unstable because both its active power and current attain astronomical values. There is also voltage, and active power dip during faults.

\begin{figure}[h!]
      \centering      \includegraphics[width=0.42\textwidth]{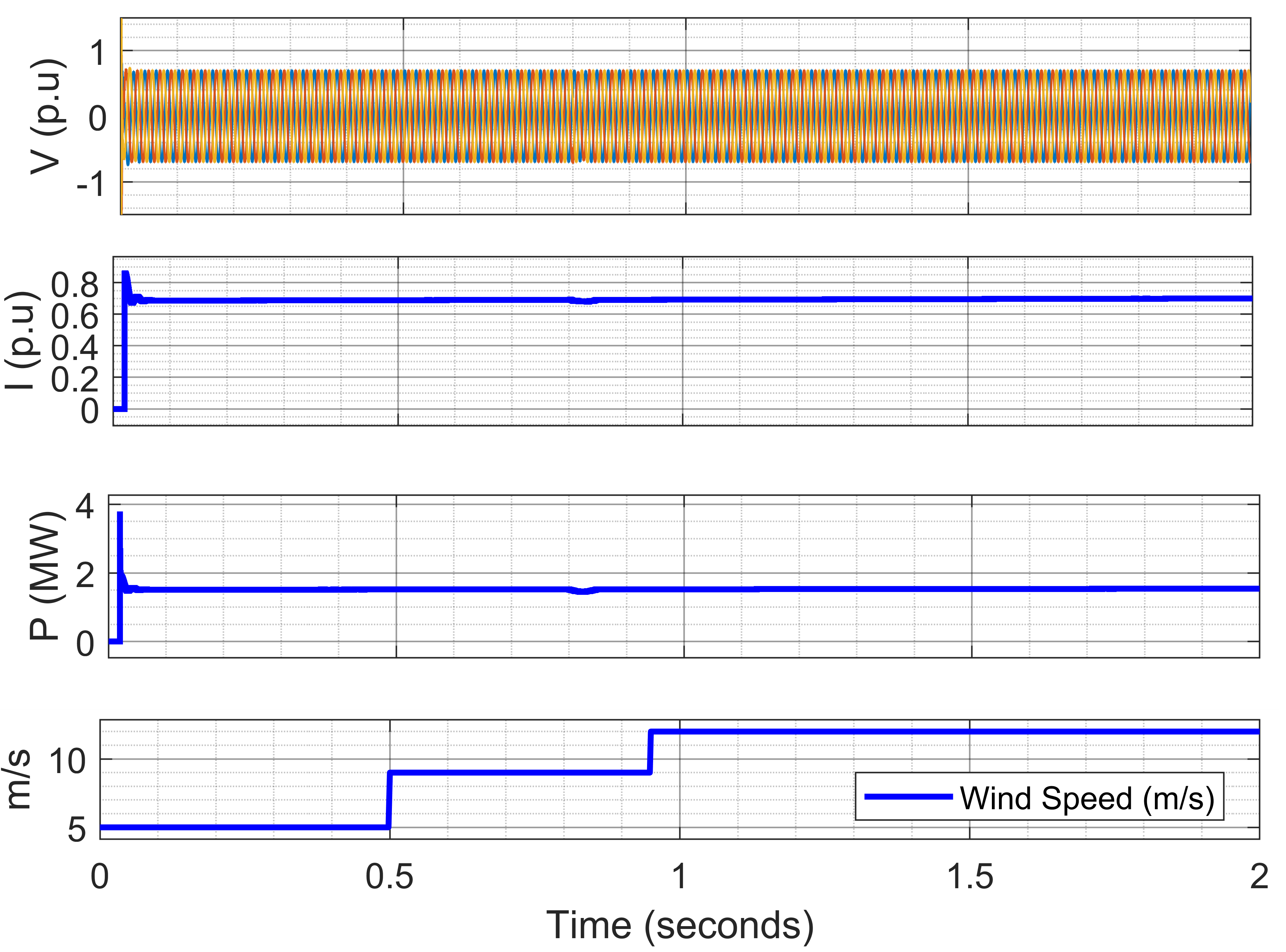}
	 \caption{WECS response with STATCOM}
     \label{stat_wecs} 
\end{figure}

\begin{figure}[t!]
      \centering
      \includegraphics[width=0.42\textwidth]{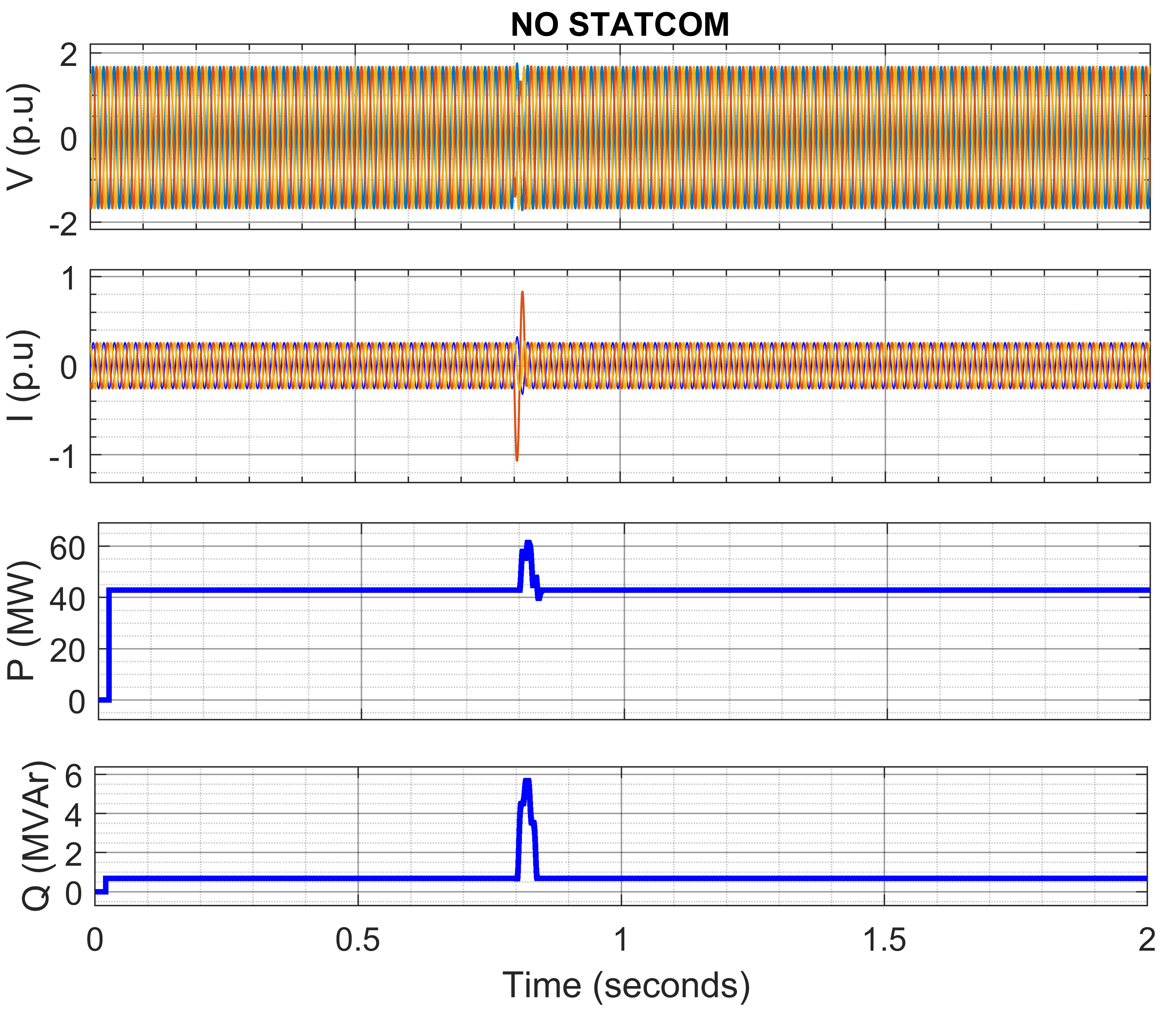}
	 \caption{PCC response with 3-phase fault without STATCOM}
     \label{No_STATCOM_bus} 
\end{figure}

Figure~\ref{No_STATCOM_bus} shows that grid voltage, current , active and reactive power become distorted between $0.8$ and $0.82\mathrm{~secs}$. Under this scenerio, there is over 30\% overvoltage because there is less reactive power and low current flow in the grid. However, the active power remains at about $42\mathrm{~MW}$ and experiences a jump only during fault, since the grid frequency is constant.

\begin{figure}[h!]
      \centering
      \includegraphics[width=0.41\textwidth]{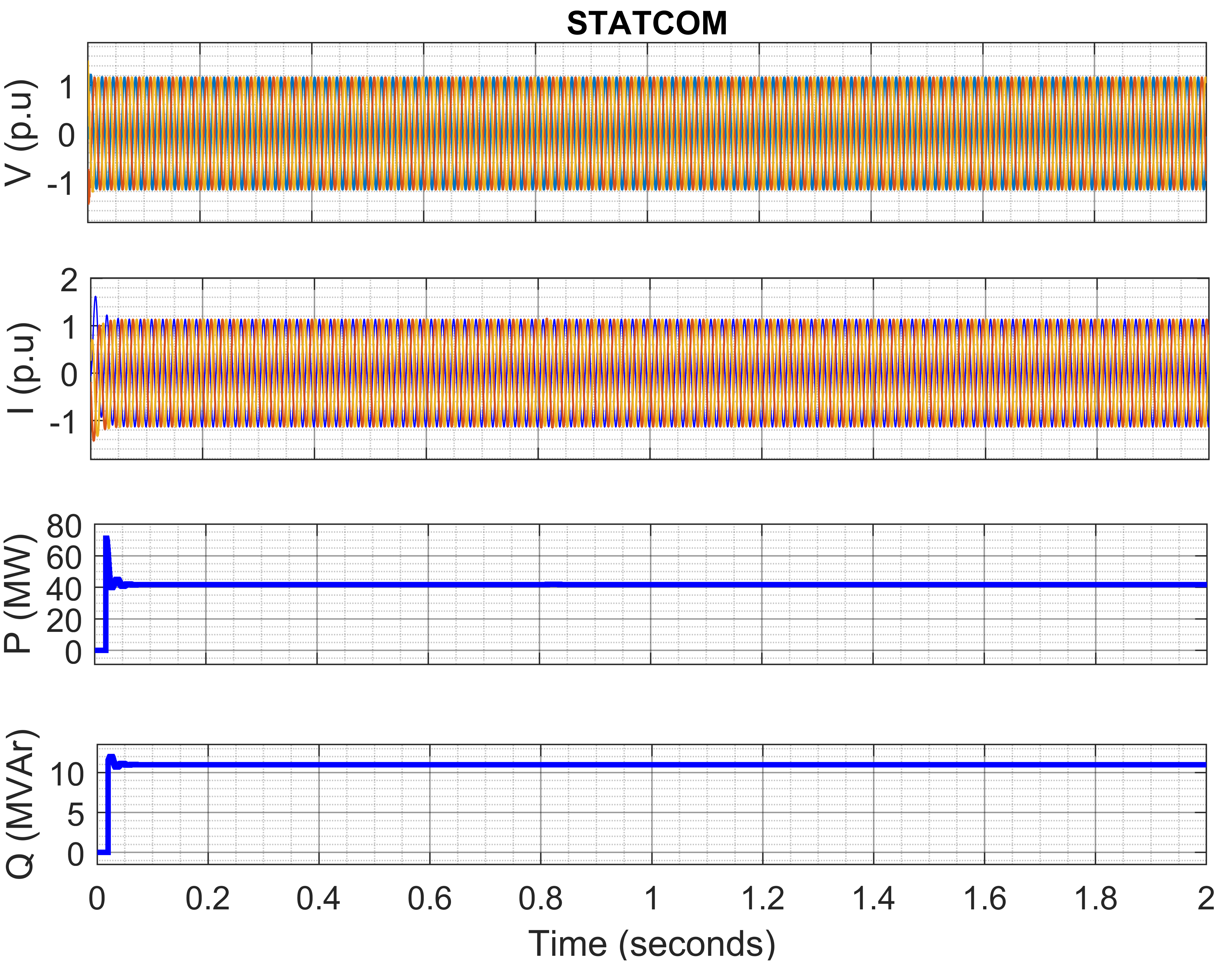}
	 \caption{PCC response with a 3-phase fault with STATCOM showing FRT capability}
     \label{STATCOM_bus} 
\end{figure}

In Figure~\ref{stat_wecs} when the STATCOM is connected to the PCC at the medium voltage level, the fault is immediately cleared and the voltage is now perfectly regulated to a nominal of $1\mathrm{~p.u}$ without no deviation. The active power of the WECS goes back to its rated value of $1.5\mathrm{~MW}$ while the current becomes regulated.
Figure~\ref{STATCOM_bus} shows that the grid current level increases along with the corresponding increase in grid reactive power to about $12\mathrm{~MVAr}$, while the grid active power remains the same as the one in the grid without STATCOM. This is simply because the grid frequency does not change. The voltage in the grid is also regulated to the nominal value of $1\mathrm{~p.u}$.

\begin{figure}[h!]
      \centering
      \includegraphics[width=0.41\textwidth]{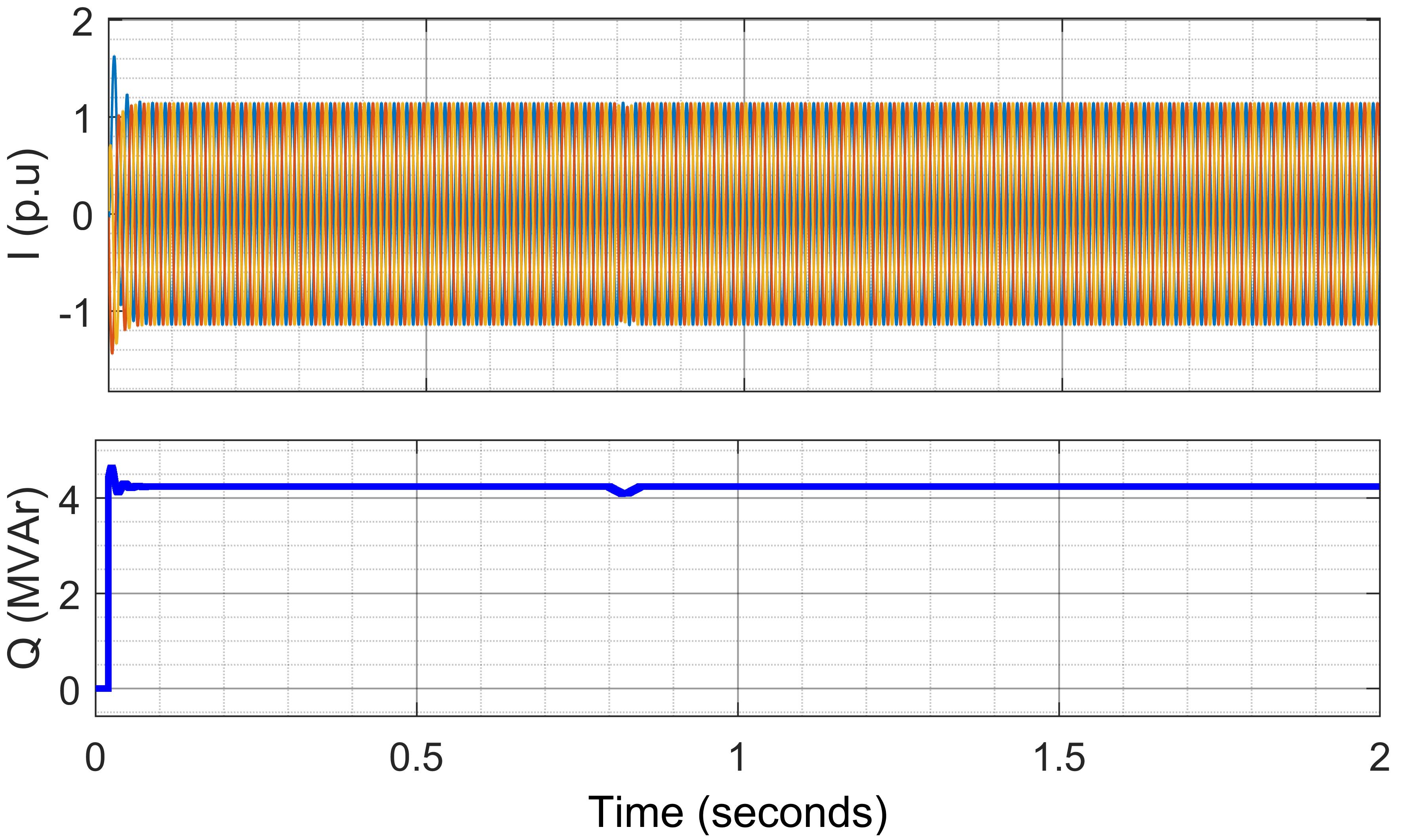}
	 \caption{Current and Reactive Power injected by STATCOM}
     \label{STATCOM_result} 
\end{figure}

Moreover, Figure~\ref{STATCOM_result} shows the current and the reactive power injected by the STATCOM to maintain the grid voltage profile and ensure fault-ride-through capabilty during faults.

\section{Conclusion}\label{sec: conclusion}
This paper illustrates the fault-ride-through capability of the grid-connected wind-energy-conversion system (WECS) equipped with STATCOM. We demonstrated that STATCOM guarantees voltage stability, reactive power compensation and robustness against injected fault. The 48 pulse configuration with zigzag connected three-phase transformers produce smooth sinusoidal signals in inverting mode. This work shows that fault may not necessarily trip WECS completely from the grid, but there may be a significant voltage deviation and the extent of response depends on the parameter reference values determined by the grid codes specific to different countries around the world. The results show that the unstable grid-tied WECS system become stable in the presence of STATCOM with the fault between $0.8$ and $0.82\mathrm{~secs}$ eliminated. Hence, with FRT capability, the WECS remain connected to the grid and also provide voltage support.
\bibliographystyle{IEEEtran}
\bibliography{references}

\begin{thebibliography}{10}
\providecommand{\url}[1]{#1}
\csname url@samestyle\endcsname
\providecommand{\newblock}{\relax}
\providecommand{\bibinfo}[2]{#2}
\providecommand{\BIBentrySTDinterwordspacing}{\spaceskip=0pt\relax}
\providecommand{\BIBentryALTinterwordstretchfactor}{4}
\providecommand{\BIBentryALTinterwordspacing}{\spaceskip=\fontdimen2\font plus
\BIBentryALTinterwordstretchfactor\fontdimen3\font minus \fontdimen4\font\relax}
\providecommand{\BIBforeignlanguage}[2]{{%
\expandafter\ifx\csname l@#1\endcsname\relax
\typeout{** WARNING: IEEEtran.bst: No hyphenation pattern has been}%
\typeout{** loaded for the language `#1'. Using the pattern for}%
\typeout{** the default language instead.}%
\else
\language=\csname l@#1\endcsname
\fi
#2}}
\providecommand{\BIBdecl}{\relax}
\BIBdecl

\bibitem{milano2018foundations}
F.~Milano, F.~D{\"o}rfler, G.~Hug, D.~J. Hill, and G.~Verbi{\v{c}}, ``Foundations and challenges of low-inertia systems,'' in \emph{2018 power systems computation conference (PSCC)}.\hskip 1em plus 0.5em minus 0.4em\relax IEEE, 2018, pp. 1--25.

\bibitem{wiser2020wind}
R.~H. Wiser, M.~Bolinger, B.~Hoen, D.~Millstein, J.~Rand, G.~L. Barbose, N.~R. Darghouth, W.~Gorman, S.~Jeong, A.~D. Mills \emph{et~al.}, ``Wind energy technology data update: 2020 edition,'' 2020.

\bibitem{operator2011wind}
A.~E.~M. Operator, ``Wind integration: International experience wp2: Review of grid codes,'' \emph{2nd October}, 2011.

\bibitem{rahimi2010grid}
M.~Rahimi and M.~Parniani, ``Grid-fault ride-through analysis and control of wind turbines with doubly fed induction generators,'' \emph{Electric Power Systems Research}, vol.~80, no.~2, pp. 184--195, 2010.

\bibitem{tarafdar2019review}
M.~Tarafdar~Hagh and T.~Khalili, ``A review of fault ride through of pv and wind renewable energies in grid codes,'' \emph{International Journal of Energy Research}, vol.~43, no.~4, pp. 1342--1356, 2019.

\bibitem{tsili2009review}
M.~Tsili and S.~Papathanassiou, ``A review of grid code technical requirements for wind farms,'' \emph{IET Renewable power generation}, vol.~3, no.~3, pp. 308--332, 2009.

\bibitem{taul2019current}
M.~G. Taul, X.~Wang, P.~Davari, and F.~Blaabjerg, ``Current reference generation based on next-generation grid code requirements of grid-tied converters during asymmetrical faults,'' \emph{IEEE Journal of Emerging and Selected Topics in Power Electronics}, vol.~8, no.~4, pp. 3784--3797, 2019.

\bibitem{rathi2005novel}
M.~R. Rathi and N.~Mohan, ``A novel robust low voltage and fault ride through for wind turbine application operating in weak grids,'' in \emph{31st Annual Conference of IEEE Industrial Electronics Society, 2005. IECON 2005.}\hskip 1em plus 0.5em minus 0.4em\relax IEEE, 2005, pp. 6--pp.

\bibitem{el2008evaluation}
H.~El-Helw and S.~B. Tennakoon, ``Evaluation of the suitability of a fixed speed wind turbine for large scale wind farms considering the new uk grid code,'' \emph{Renewable Energy}, vol.~33, no.~1, pp. 1--12, 2008.

\bibitem{saad1998application}
Z.~Saad-Saoud, M.~Lisboa, J.~Ekanayake, N.~Jenkins, and G.~Strbac, ``Application of statcoms to wind farms,'' \emph{IEE Proceedings-Generation, Transmission and Distribution}, vol. 145, no.~5, pp. 511--516, 1998.

\bibitem{muyeen2005stabilization}
S.~Muyeen, M.~A. Mannan, M.~H. Ali, R.~Takahashi, T.~Murata, and J.~Tamura, ``Stabilization of grid connected wind generator by statcom,'' in \emph{2005 International Conference on Power Electronics and Drives Systems}, vol.~2.\hskip 1em plus 0.5em minus 0.4em\relax IEEE, 2005, pp. 1584--1589.

\bibitem{kuiava2009control}
R.~Kuiava, R.~A. Ramos, and N.~G. Bretas, ``Control design of a statcom with energy storage system for stability and power quality improvements,'' in \emph{2009 IEEE International Conference on Industrial Technology}.\hskip 1em plus 0.5em minus 0.4em\relax IEEE, 2009, pp. 1--6.

\bibitem{teleke2011application}
S.~Teleke, A.~Yazdani, B.~Gudimetla, J.~Enslin, and J.~Castaneda, ``Application of statcom for power quality improvement,'' in \emph{2011 IEEE/PES Power Systems Conference and Exposition}.\hskip 1em plus 0.5em minus 0.4em\relax IEEE, 2011, pp. 1--6.

\bibitem{barrado2010power}
J.~Barrado, R.~Grino, and H.~Valderrama-Blavi, ``Power-quality improvement of a stand-alone induction generator using a statcom with battery energy storage system,'' \emph{IEEE transactions on power delivery}, vol.~25, no.~4, pp. 2734--2741, 2010.

\bibitem{hingorani1999understanding}
N.~G. Hingorani and L.~Gyugyi, ``Understanding facts,'' \emph{(No Title)}, 1999.

\bibitem{adetokun2021application}
B.~B. Adetokun and C.~M. Muriithi, ``Application and control of flexible alternating current transmission system devices for voltage stability enhancement of renewable-integrated power grid: A comprehensive review,'' \emph{Heliyon}, vol.~7, no.~3, 2021.

\bibitem{ledesma2005doubly}
P.~Ledesma and J.~Usaola, ``Doubly fed induction generator model for transient stability analysis,'' \emph{IEEE transactions on energy conversion}, vol.~20, no.~2, pp. 388--397, 2005.

\bibitem{nemmour2010advanced}
A.~Nemmour, F.~Mehazzem, A.~Khezzar, M.~Hacil, L.~Louze, and R.~Abdessemed, ``Advanced backstepping controller for induction generator using multi-scalar machine model for wind power purposes,'' \emph{Renewable Energy}, vol.~35, no.~10, pp. 2375--2380, 2010.

\bibitem{sakamoto2007output}
R.~Sakamoto, T.~Senjyu, T.~Kaneko, N.~Urasaki, T.~Takagi, and S.~Sugimoto, ``Output power leveling of wind turbine generator by pitch angle control using h-infinity control,'' \emph{IEEJ Transactions on Power and Energy}, vol. 127, no.~1, pp. 86--93, 2007.

\bibitem{jurado2003adaptive}
F.~Jurado and J.~R. Saenz, ``An adaptive control scheme for biomass-based diesel--wind system,'' \emph{Renewable energy}, vol.~28, no.~1, pp. 45--57, 2003.

\bibitem{singh2012new}
B.~Singh and V.~S. Kadagala, ``A new configuration of two-level 48-pulse vscs based statcom for voltage regulation,'' \emph{Electric Power Systems Research}, vol.~82, no.~1, pp. 11--17, 2012.

\bibitem{abdellaoui2015statcom}
A.~Abdellaoui, A.~Yangui, A.~Saidi, and H.~H. Abdallah, ``Statcom-based 48-pulses three level gto dedicated to var compensation and power quality improvement,'' in \emph{2015 International Conference on Sustainable Mobility Applications, Renewables and Technology (SMART)}.\hskip 1em plus 0.5em minus 0.4em\relax IEEE, 2015, pp. 1--7.

\bibitem{park1929two}
R.~H. Park, ``Two-reaction theory of synchronous machines generalized method of analysis-part i,'' \emph{Transactions of the American Institute of Electrical Engineers}, vol.~48, no.~3, pp. 716--727, 1929.

\bibitem{miller2003dynamic}
N.~W. Miller, W.~W. Price, and J.~J. Sanchez-Gasca, ``Dynamic modeling of ge 1.5 and 3.6 wind turbine-generators,'' \emph{GE-Power systems energy consulting}, no. 3.0, 2003.

\bibitem{ameli2022hierarchical}
S.~Ameli and O.~M. Anubi, ``Hierarchical robust control for variable-pitch wind turbine with actuator faults,'' \emph{International Journal of Robust and Nonlinear Control}, vol.~32, no.~12, pp. 7039--7056, 2022.

\bibitem{10252896}
A.~Olajube and O.~M. Anubi, ``Model-based loss minimization control of a squirrel cage induction motor drive with shorted rotor under indirect field orientation,'' in \emph{2023 IEEE Conference on Control Technology and Applications (CCTA)}, 2023, pp. 759--765.

\bibitem{sahoo2006modeling}
A.~K. Sahoo, K.~Murugesan, and T.~Thygarajan, ``Modeling and simulation of 48-pulse vsc based statcom using simulink’s power system blockset,'' in \emph{2006 India International Conference on Power Electronics}.\hskip 1em plus 0.5em minus 0.4em\relax IEEE, 2006, pp. 303--308.

\bibitem{suresh2010dynamic}
Y.~Suresh and A.~Panda, ``Dynamic performance of statcom under line to ground faults in power system,'' 2010.

\end{thebibliography}

\end{document}